\documentclass[11pt]{amsart}
\usepackage{amsmath, amsfonts, amsthm, amssymb, multicol, mathtools, dsfont, verbatim}
\usepackage{graphicx}
\usepackage{float, hyperref}
\usepackage{scalerel,stackengine, subcaption}
\usepackage[usenames,dvipsnames]{xcolor}
\usepackage{enumitem}
\usepackage{pgfplots}
\usepgfplotslibrary{fillbetween}
\pgfplotsset{width=10cm,compat=1.9}
\usepackage[square,sort,comma,numbers]{natbib}
\allowdisplaybreaks

\usepackage[square,sort,comma,numbers]{natbib}
\usepackage{fancyhdr}

\makeatletter
\def\@setauthors{%
  \begingroup
  \def\thanks{\protect\thanks@warning}%
  \trivlist
  \centering\footnotesize \@topsep30\p@\relax
  \advance\@topsep by -\baselineskip
  \item\relax
  \author@andify\authors
  \def\\{\protect\linebreak}

  \normalsize\lowercase{\authors}%
  
	\ifx\@empty\contribs
  \else
    ,\penalty-3 \space \@setcontribs
    \@closetoccontribs
  \fi
  \endtrivlist
  \endgroup
}
\def\@settitle{\begin{center}
\LARGE\lowercase{\@title}
  \end{center}%
}
\makeatother

\definecolor{lightblue}{HTML}{2B77A4}
\definecolor{darkred}{HTML}{9E0D0D}
\hypersetup{
	colorlinks=true,
	linkcolor=darkred,
	urlcolor=darkred,
	citecolor=lightblue
}
\numberwithin{equation}{section}

\newtheorem{thm}{Theorem}[section]

\newtheorem{cor}[thm]{Corollary}

\newtheorem{prop}[thm]{Proposition}

\newtheorem*{thmA}{Theorem A}

\renewcommand{\epsilon}{\varepsilon}

\newcommand{\rd}{\mathbb{R}^d}
\newcommand{\ff}{\mathbb{F}_q}

\renewcommand{\geq}{\geqslant}
\renewcommand{\leq}{\leqslant}

\title{Fourier analytic properties of  Kakeya sets in finite fields}

\author{Jonathan M. Fraser\\ \\
 University of St Andrews, Scotland\\
\MakeLowercase{Email: jmf32@st-andrews.ac.uk}}
\thanks{The  author was    financially supported by a  \emph{Leverhulme Trust Research Project Grant} (RPG-2023-281),  an \emph{EPSRC Standard Grant} (EP/Y029550/1), and an \emph{EPSRC Open Fellowship} (EP/Z533440/1).}

\begin{document}


\maketitle
\thispagestyle{empty}

\begin{abstract}
We prove that a Kakeya set in a vector space over a  finite field of size $q$ always supports a probability measure whose Fourier transform is bounded by $q^{-1}$ for all non-zero frequencies.  We show that this  bound is sharp in all dimensions at least 2.  In particular,  this provides a new and self-contained proof that a Kakeya set in dimension 2 has size at least $q^2/2$    (which is  asymptotically sharp).  We also establish analogous results for sets containing $k$-planes in a given set of  orientations.
\end{abstract}

\tableofcontents

\section{Introduction}

 \subsection{Kakeya sets in finite fields}

A  \textit{Kakeya set} in $\mathbb{R}^d$ is a set which contains a unit line segment in every possible direction.  The   \emph{Kakeya conjecture} is that such sets must have Hausdorff dimension $d$. This is despite the fact that they can have zero $d$-dimensional Lebesgue measure for $d \geq 2$; a result of Besicovitch \cite{besicovitch}.  Davies solved the planar ($d=2$) case  \cite{davies} and the $d=3$ case was recently resolved in a breakthrough paper of Wang and Zahl \cite{wangzahl}.    The problem   is open for $d \geq 4$.

  Oberlin proved that Kakeya sets in $\mathbb{R}^2$ have \emph{Fourier} dimension 2, that is, Kakeya sets in $\mathbb{R}^2$  are Salem sets \cite{oberlin}.  Oberlin's result implies (and is rather stronger than) Davies' result that such Kakeya sets have \emph{Hausdorff} dimension 2. We generalised Oberlin's argument in \cite{cones} to prove that the Fourier dimension of a Kakeya set in $\rd$ is \emph{at least} 2   and in \cite{product} we proved that there exist Kakeya sets in $\rd$ with Fourier dimension precisely 2 for all $d \geq 3$.  In particular, the fact that Kakeya sets in the plane are Salem does not generalise to higher dimensions.

  Wolff  proposed the following finite field analogue of the Kakeya problem \cite{wolffprop}. Let $\ff$ be the finite field of $q$ elements where $q$ is a power of a prime, and let $\ff^d$ be the $d$-dimensional vector space over $\ff$, where $d \geq 2$ is an integer. In this setting a set $K \subseteq \mathbb{F}_q^d$ is called a \emph{Kakeya set} if it contains a line in every direction, that is, for all $x \in \mathbb{F}_q^d$, there exists $y \in \mathbb{F}_q^d$ such that
\[
\{ y+a \,  x : a \in \mathbb{F}_q \} \subseteq K.
\]
The \emph{finite fields Kakeya conjecture} is  that   Kakeya sets $K \subseteq \mathbb{F}_q^d$, must have cardinality at least a constant multiple of $q^d$.  This conjecture  was sensationally solved by Dvir \cite{dvir} in full generality, with the case $d=2$ having been settled  by Wolff \cite{wolffprop}.  

To the best of our knowledge the finite field analogue of Oberlin's result is unknown.  We prove it here and in the process give a new proof of the asymptotically  sharp lower bound  for the size of a Kakeya set  in $\mathbb{F}_q^2$.
\begin{thmA}
If $K \subseteq \mathbb{F}_q^2$ is a   Kakeya set, then it is a Salem set of size at least  $q^2/2$.
\end{thmA}
 This result is a special case of our main results which are discussed in Section \ref{results}; see Theorem \ref{main} and  Corollary \ref{maind=2}.  Also, see Section \ref{dfa} for the definition of a Salem set in the finite fields setting.
 
 A natural generalisation of a Kakeya set is a $(d,k)$-set, which is a set in $d$-dimensional space  which contains a $k$-plane in every orientation.  There has been a lot of work on these sets in the Euclidean setting, but we only discuss the finite fields setting here, see \cite{ellenberg, kopparty}.  We consider a slight generalisation of $(d,k)$-sets where we only require a given subset of possible orientations.  Write $G(d,k)$ for the set of all $k$-dimensional subspaces of $\ff^d$ and let $\Gamma \subseteq G(d,k)$.  We say $K \subseteq \ff^d$ is a $(d,k, \Gamma)$\emph{-set} if for all $V \in \Gamma$, there exists $u \in \ff^d$ such that $u+V \subseteq K$.  In particular, $K$ is a  $(d,1, G(d,1))$-set if and only if it is a Kakeya set.   We first obtain Fourier analytic estimates for $(d,k, \Gamma)$-sets in Theorem \ref{maindk} before specialising to $(d,k)$-sets in Corollary \ref{maindkcor} and then Kakeya sets in Proposition \ref{mainex}, Theorem \ref{main}, and Corollary \ref{maind=2}.

 \subsection{Notation}

Throughout the paper we write $A \lesssim B$ to mean there is a constant $c$  independent of $q$ for which $A \leq cB$.  Similarly we write $A \gtrsim B$ to mean $B \lesssim A$ and $A \approx B$ if $A\lesssim B$ and $A \gtrsim B$.  One should always think of $q$ as being very large compared to any implicit constants. We also write $| E|$ to denote the cardinality of a set $E$ and $A \sim B$ to mean that $A/B \to 1$ as $q \to \infty$.


\subsection{Fourier analysis in finite fields} \label{dfa}

 The \emph{Fourier transform}   of a function $f : \ff^d \to \mathbb{C}$ is the function  $\widehat f :  \ff^d \to \mathbb{C}$   defined by
\[
\widehat f(\xi) \coloneqq   \sum_{x \in \ff^d} \chi(-\xi  \cdot x)f(x),
\]
where $\chi : \ff \to S^1 \subseteq \mathbb{C}$ is a non-trivial additive character, which we fix from now on. The specific choice of $\chi$ does not play an important role in what follows. Above, $\xi \cdot x$ denotes the standard dot product on $\ff^d$, which takes values in $\ff$.

 \emph{Plancherel's formula} gives a useful $L^2$ identity connecting a function and its Fourier transform.  It states that, for  all $f : \ff^d \to \mathbb{C}$,
\begin{equation}
\label{plancherel}
 \sum_{\xi \in \ff^d} |\widehat f(\xi)|^2  = q^d \sum_{x \in \ff^d}  |f(x)|^2.
\end{equation}
 A probability measure $\mu$ on $\ff^d$ is simply a non-negative real-valued function that sums to 1.  Let $E\subseteq \ff^d$ be the support of a probability measure $\mu$, that is, $E= \{ x \in \ff^d : \mu(x) >0\}$.  Then, for all frequencies $\xi \in \ff^d$,
 \[
 |\widehat \mu (\xi)| \leq  |\widehat \mu (0)| = 1
 \]
 and, by Plancherel's formula \eqref{plancherel},
\begin{equation} \label{simple1}
 \sum_{\xi \in \ff^d}  |\widehat \mu (\xi)|^2 = q^d  \sum_{x\in \ff^d}  | \mu (x)|^2 \geq q^d \sum_{x\in E}  |E|^{-2} = q^d |E|^{-1}
  \end{equation}
 where we used here the simple fact that the $L^2$ norm of a finitely supported probability distribution is minimised by the uniform distribution.  Indeed, if a probability measure on a finite set $E$ is not uniform, then there must exist weights $a,b>0$ with  $ a<|E|^{-1} <b$.  But then the $L^2$ norm can be reduced by replacing $a$ and $b$ with $|E|^{-1}$ and $a+b-|E|^{-1}$ since
 \[
 a^2+b^2 > |E|^{-2}+(a+b-|E|^{-1})^2 \quad \Leftrightarrow \quad  (|E|^{-1} -a)(|E|^{-1}-b) <0
 \]
 and finitely many iterations of this procedure yields the uniform distribution as the unique minimiser.
   We also have the simple estimate
\begin{equation} \label{simple2}
  \sum_{\xi \in \ff^d}  |\widehat \mu (\xi)|^2  \leq 1 + (q^d-1)\sup_{\xi \in \ff^d \setminus \{0\}} |\widehat \mu (\xi)|^2.
  \end{equation}
It follows from \eqref{simple1} and \eqref{simple2} that
\begin{equation} \label{lowerb1}
 \sup_{\xi \in \ff^d \setminus \{0\}} |\widehat \mu (\xi)|  \geq \sqrt{\frac{q^d |E|^{-1}-1}{q^d-1}}.
 \end{equation}
In particular, if  $|E| \leq cq^d$ for some $c \in (0,1)$ (that is, if $|\ff^d \setminus E| \gtrsim q^d$), then
\begin{equation} \label{lowerb}
 \sup_{\xi \in \ff^d \setminus \{0\}} |\widehat \mu (\xi)|  \geq  
  \sqrt{1-c} \,  |E|^{-1/2} \approx |E|^{-1/2}.
 \end{equation}
This allows us to get cardinality estimates for a set via the Fourier transform of measures it supports.  Indeed, if there exists a probability measure supported on a set $E$ such that 
\[
 \sup_{\xi \in \ff^d \setminus \{0\}} |\widehat \mu (\xi)|  \leq C q^{-\beta/2}
\]
 for some $C \geq 0$ and $\beta \in (0,d]$, then \eqref{lowerb1} gives
\begin{equation} \label{sizeest}
 |E| \geq \frac{q^d}{C^2q^{-\beta}(q^d-1)+1} \gtrsim q^\beta
\end{equation}
 Further, if $\beta>d$ then one may conclude $|E| \sim q^d$ and if $\beta>2d$ then $E = \ff^d$ for large enough $q$. 
 
Given the above observations, we say $E \subseteq \ff^d$ is \emph{Salem} if there exists a probability measure $\mu$ with support contained in  $E$ such that 
 \[
 \sup_{\xi \in \ff^d \setminus \{0\}} |\widehat \mu (\xi)|  \lesssim |E|^{-1/2}.
 \]
 We note  that this is \emph{not} the same as the definition of a Salem set in $\ff^d$ introduced by Iosevich--Rudnev \cite{iosevich}. In \cite{iosevich} the measure $\mu$ was required to be the `surface measure', or uniform distribution, on $E$.  We suggest that the definition of Salem in \cite{iosevich} should be thought of as \emph{strongly Salem} or perhaps \emph{uniformly Salem}. In particular, in Euclidean space one is free to demonstrate that a given set is Salem by choosing an arbitrary Borel measure supported on the set and therefore  our definition above is a more appropriate analogue of the definition in Euclidean space, see \cite{mattila}.  That said, certain uniformity properties are forced upon a measure $\mu$ witnessing that a set $E$ (with $|E| \leq c q^d$ for some $c \in (0,1)$) is Salem.  In particular, $E$ must be supported on a subset of $E$ with cardinality $\approx |E|$ and,
 by virtue of \eqref{simple1}, we must have
 \[
  \sum_{x\in \ff^d}  | \mu (x)|^2 \approx \sum_{x\in E}  |E|^{-2} = |E|^{-1} 
 \]
 and so perhaps the definitions are not as different as they first appear (although they are clearly not equivalent).

\section{Main results} \label{results}

Our first result proves that  $(d,k, \Gamma)$-sets always support  a probability measure with polynomial Fourier decay, provided $|\Gamma| \gtrsim q^{k(d-k-1)+\varepsilon}$ for some $\varepsilon>0$.  It is useful to keep in mind that
\[
q^{k(d-k)} \leq |  G(d,k)| \leq  q^{k(d-k)}(1+o(1))
\]
and
\[
\frac{|  G(d-1,k)|}{|  G(d,k)|} \leq q^{-k}.
\]
Both of these facts can be derived easily by considering Gaussian binomial coefficients; we omit the details.

\begin{thm} \label{maindk}
Let   $d>k \geq 1$ be integers and $\Gamma \subseteq G(d,k)$.  Suppose $K \subseteq \ff^d$ is such that for all $V \in \Gamma$, there exists $u \in \ff^d$ such that $u+V \subseteq K$.  Then there exists a probability measure $\mu$ on $K$ satisfying
\[
|\widehat \mu (\xi)| \leq |\Gamma|^{-1} |  G(d-1,k)| \leq q^{-k} \, \frac{|  G(d,k)| }{|\Gamma|} 
\]  
for all non-zero $\xi \in \ff^d$.   In particular, $|K| \gtrsim \min\{q^d, |\Gamma|^2q^{-2k(d-k-1)}\}$ and more precise estimates are available from \eqref{sizeest} if required.
\end{thm}

We delay the proof of Theorem \ref{maindk} until Section \ref{proof1}. Theorem \ref{maindk} has no content when $|\Gamma| \lesssim q^{k(d-k-1)}$ but this is for good reason. In particular,   $V_0 \in G(d,d-1)$ is a $(d,k, \Gamma_0)$-set where $\Gamma_0 = \{ V \in G(d,k) : V \subseteq V_0\}$.  Then $|\Gamma| \approx q^{k(d-k-1)}$ but $V_0$ cannot support a probability measure with non-trivial bounds for the Fourier transform at $\xi \in V_0^\perp$.  

Setting $\Gamma=G(d,k)$  we get the following corollary concerning traditional $(d,k)$-sets. 

\begin{cor} \label{maindkcor}
Let   $d>k \geq 1$ be integers and   $K \subseteq \ff^d$ be a $(d,k)$-set.  Then there exists a probability measure $\mu$ on $K$ satisfying
\[
|\widehat \mu (\xi)| \leq q^{-k} 
\]  
for all non-zero $\xi \in \ff^d$.  In particular, $|K| \gtrsim \min\{q^d, q^{2k}\}$ and if $k = d/2$, then $|K| \approx  q^d$, and if $k>d/2$, then $|K| \sim q^d$.
\end{cor}
The conclusion that $(d,k)$-sets asymptotically take up the whole space when $k > d/2$ is a finite field analogue of a result of Falconer in the Euclidean setting \cite{falconer}, which states that $(d,k)$-sets in $\mathbb{R}^d$ must have positive Lebesgue measure when $k>d/2$.  In the finite fields setting, the cardinality estimate in Corollary \ref{maindkcor} is known, and in fact much stronger estimates are available.  Therefore, the main novelty of Corollary \ref{maindkcor} is in the Fourier analytic information, but we do  offer a new proof of the optimal asymptotic cardinality estimate in the case $k > d/2$.  In fact, for all $k \geq 2$, it is known that  $(d,k)$-sets in $\ff^d$ have cardinality $\sim q^d$, see \cite{ellenberg, kopparty}.

We now focus on the case $k=1$, that is, the setting of Kakeya sets.  We first prove that the estimates above are sharp in all ambient dimensions $d \geq 2$.  

\begin{prop} \label{mainex}
For all   $d \geq 2$,  $\kappa \in (0,1)$, and sufficiently large $q$,  there exists a Kakeya set $K \subseteq \ff^d$ such that for all probability measures $\mu$ on $K $  
\[
|\widehat \mu (\xi)| \geq \kappa \,  q^{-1}.
\]  
\end{prop}

We delay the proof of Proposition  \ref{mainex} until Section \ref{proof2}. Combining Theorem \ref{maindk} and Proposition \ref{mainex}, we get the following theorem, which we see as our main result.  In particular, in the case $d=2$ this gives a finite fields analogue of Oberlin's result that Kakeya sets in the Euclidean plane are Salem \cite{oberlin},  a  strengthening of Dvir's result \cite{dvir}, and  a new self-contained proof that Kakeya sets in $\ff^2$ must have cardinality $\approx q^2$. In higher dimensions, we do not obtain a proof of the finite fields Kakeya conjecture, and Proposition \ref{mainex} shows that it cannot be proved this way for $d>2$.

\begin{thm} \label{main}
Let $d \geq 2$ and $K \subseteq \ff^d$ be a Kakeya set.  Then there exists a probability measure $\mu$ on $K$ satisfying
\[
|\widehat \mu (\xi)| \leq  q^{-1}
\]  
for all non-zero $\xi \in \ff^d$.  In particular, $|K| \gtrsim q^2$.  Moreover, for all $d \geq 2$,  $\kappa \in (0,1)$, and sufficiently large $q$,  there exists a Kakeya set $K' \subseteq \ff^d$ such that for all probability measures $\mu$ on $K'$ 
\[
|\widehat \mu (\xi)| \geq \kappa \,  q^{-1}
\]  
and so the Fourier analytic bound above is sharp for all $d \geq 2$.
\end{thm}

The special case $d=2$ is worth further scrutiny since there we get the right cardinality estimate.  In fact,  we obtain a new proof of the sharp density result for Kakeya sets in dimension 2. This follows from the Fourier bound in Theorem \ref{main} and the more precise cardinality estimate \eqref{sizeest}.

\begin{cor} \label{maind=2}
If   $K \subseteq \ff^2$ is a Kakeya set, then
\[
|K| \geq  \frac{q^2}{2-q^{-2}} \geq q^2/2.
\]
\end{cor}

Corollary \ref{maind=2} improves on the estimate $|K| \geq q^2/4$ originally proved by Wolff \cite{wolffprop} and also again as the $d=2$ case of the estimate  $|K| \geq 2^{-d}q^d$ for Kakeya sets in $\ff^d$  from \cite{smallkakeya}.  It also improves upon (but is asymptotically equivalent to) the bound $|K| \geq q(q-1)/2$ obtained by Dvir \cite{dvir}.  The state-of-the-art   was recently improved to
\begin{equation} \label{bukhbound}
|K| \geq  \frac{q^d}{(2-q^{-1})^{d-1}} \geq 2^{1-d}q^d  
\end{equation}
 for Kakeya sets in $\ff^d$ in  \cite{bukh} and this result served to establish the sharp density result for Kakeya sets in dimensions $d \geq 3$ for the first time; the $d=2$ case having been settled by Dvir.   The fact that density $2^{1-d}$ cannot be beaten   was  proved  by Dvir who constructed  a Kakeya set $K \subseteq \ff^d$ with $|K| \leq 2^{1-d}q^d  +O(q^{d-1})$, see \cite{saraf}. In particular, Corollary \ref{maind=2} provides a novel proof of the asymptotically sharp (density) result for  Kakeya sets in $\ff^2$, although the bound we get  is slightly weaker than (although asymptotically equivalent to)    \eqref{bukhbound}.

\section{Proofs}

\subsection{Proof of Theorem \ref{maindk}} \label{proof1}

Let $K \subseteq \ff^d$ be a $(d,k, \Gamma)$-set.   For each $V \in \Gamma$, fix a minimal basis $\{v_1^V, \dots, v_k^V\} \subseteq V$.  Also, write
\[
V^\perp = \{ u \in \ff^d : u \cdot v = 0 \text{ for all $v \in V$}\}.
\]
By definition, for some $\{u_V  : V \in \Gamma\} \subseteq \ff^d$ the set 
\[
K' = \bigcup_{V \in \Gamma} \bigcup_{a_1, \dots, a_k \in \mathbb{F}_q} \left\{u_V+\sum_{i=1}^ka_i v_i^V\right\}
\]
is contained in $K$.  Define $\phi: \ff^d \to \mathbb{Z}$ by 
\[
\phi(z) = \left\lvert \left\{ V \in  \Gamma : z \in  u_V+\sum_{i=1}^ka_i v_i^V \right\}\right\rvert,
\]
that is, $\phi(z)$ is the number of incidences between the affine $k$-planes defining $K'$ and the point $z$.  Clearly $\phi(z) = 0$ if and only if  $z \notin K'$ and therefore setting
\[
\mu(z) = 
\frac{\phi(z)}{\sum_{u \in K'} \phi(u)}  
\]
defines a probability measure $\mu$ on $K$ (with support equal to $K' \subseteq K$).   Moreover,
\begin{equation} \label{weight}
\sum_{u \in K'} \phi(u)  = |  \Gamma| \,  q^k   .
\end{equation}
  Writing $\mathbf{1}$ for the indicator function on $\ff$, note that
\begin{equation} \label{indicator}
\widehat{\mathbf{1} } (u) = \left\{ \begin{array}{cc}
0 & u \neq 0\\
q & u = 0.
\end{array}\right.  
\end{equation}
Also, note that for all non-zero $\xi \in \ff^d$,
\begin{equation} \label{complement}
|\{V \in \Gamma: \xi \in V^\perp\} | \leq |\{V \in G(d,k): \xi \in V^\perp\} | = |  G(d-1,k)|.
\end{equation}
Then
\begin{align*}
|\widehat \mu(\xi)| 
& = q^{-k } |\Gamma|^{-1} \left\lvert   \sum_{z \in \mathbb{F}_q^d} \chi(-\xi \cdot z) \phi(z)\right\rvert \qquad \text{(by \eqref{weight})}\\
& =  q^{-k } |\Gamma|^{-1} \left\lvert  \sum_{V \in \Gamma} \sum_{a_1, \dots, a_k \in\mathbb{F}_q}  \chi(-\xi \cdot u_V) \chi\left(-\xi \cdot  \left(\sum_{i=1}^ka_i v_i^V\right)\right)   \right\rvert \qquad \text{(by definition of $\phi$)}\\
&  \leq  q^{-k } |\Gamma|^{-1}  \sum_{V \in \Gamma}\left\lvert  \sum_{a_1, \dots, a_k \in\mathbb{F}_q}  \chi\left(-\sum_{i=1}^k\xi \cdot(a_i  v_i^V)\right) \right\rvert\\
&  \leq  q^{-k } |\Gamma|^{-1} \sum_{V \in \Gamma} \sum_{a_1, \dots, a_k \in\mathbb{F}_q} \prod_{i=1}^k\left\lvert  \chi\left(- a_i (\xi \cdot v_i^V)\right) \right\rvert\\
& =  q^{-k } |\Gamma|^{-1}  \sum_{V \in \Gamma}   \prod_{i=1}^k | \widehat{\mathbf{1}}(\xi \cdot v_i^V) | \\
& =   q^{-k } |\Gamma|^{-1}  \cdot  |\{V \in \Gamma: \xi \in V^\perp\} |  \cdot q^k \qquad \text{(by \eqref{indicator})}\\
&\leq   |\Gamma|^{-1}   |  G(d-1,k)| \qquad \text{(by \eqref{complement})} 
\end{align*}
as required.    Finally, it follows that $|K| \geq |K'| \gtrsim \min\{q^d, |\Gamma|^2q^{-2k(d-k-1)}\}$ by applying \eqref{sizeest}.

\subsection{Proof of Proposition \ref{mainex}} \label{proof2}

Let $c \in (1/2,1)$ and let $K_0 \subseteq \mathbb{F}_q^2$ be a   Kakeya set with $|K_0| \leq c q^2$ for large enough $q$.  Such Kakeya sets exist by \cite{saraf};  recall also Corollary \ref{maind=2} and \cite{bukh}.  Then $K = K_0 \oplus \ff^{d-2} $ is     a  Kakeya set in $\ff^d \equiv\ff^2 \oplus \ff^{d-2} $.   Let $\mu$ be an arbitrary probability measure on $K$ and let $\nu$ be the projection of $\mu$ into $\ff^2$ defined by
\[
\nu(x) = \sum_{y \in  \ff^{d-2}}  \mu(x,y).
\]
Then $\nu$ is a probability measure supported on $K_0$. Let $\xi=(\xi_1, 0)$ where $\xi_1 \in \ff^2$ and $0 \in \ff^{d-2}$.  Then
\[
\widehat \mu (\xi) =  \sum_{(x,y) \in \ff^2 \oplus \ff^{d-2} \equiv \mathbb{F}_q^d} \chi(-\xi_1 \cdot x) \mu(x,y) = \sum_{ x  \in \ff^2}   \chi(-\xi_1 \cdot x) \nu(x) = \widehat \nu(\xi_1).
\]
Since $|K_0| \leq c q^2$   it follows from \eqref{lowerb} that  
\[
\sup_{\xi_1 \in \ff^2 \setminus \{0\}} |\widehat \nu(\xi_1)| \geq \sqrt{1-c} \, |K_0|^{-1/2} \geq  \sqrt{\frac{1-c}{c}} \, q^{-1}
\]
and therefore  
\[
\sup_{\xi  \in \ff^d \setminus \{0\}} |\widehat \mu(\xi)| \geq  \sqrt{\frac{1-c}{c}} \, q^{-1}
\]
and the result follows since the constant can be made arbitrarily close to 1 by letting $c \searrow 1/2$.

\end{document}